\documentclass[12pt,a4paper,oneside]{article}
\usepackage{amsmath}
\usepackage[latin1]{inputenc}
\usepackage{amsmath}
\usepackage{amsfonts}
\usepackage{amssymb}
\usepackage{latexsym}

\newtheorem{theorem}{Theorem}[section]
\newtheorem{proposition}[theorem]{Proposição}
\newtheorem{lemma}[theorem]{Lema}

\def \fimprova {\blacksquare}

\begin{document}
\title{\textbf{Finitely Generated Groups $G$ such that $G/Z(G) \simeq C_p \times C_p$}}
\author{Mariana Garabini Cornelissen \footnote{Corresponding author - mariana@ufsj.edu.br} \\ \small\emph{{Universidade Federal de São João Del Rei, Brazil}}
\\ \\ C.Polcino Milies \\ \small\emph{{Universidade de São Paulo, Brazil}}}
\date{}
\maketitle


\begin{abstract}
Finite groups $G$ such that $G/Z(G) \simeq C_2 \times C_2$ where $C_2$ denotes a cyclic group of order 2
and $Z(G)$ is the center of $G$ were studied in \cite{casofinito} and were used to classify finite loops with alternative loop algebras. In this paper we extend this result to finitely generated groups such that $G/Z(G) \simeq C_p \times C_p$ where $C_p$ denotes a cyclic group of prime order $p$ and provide an explicit description of all such groups.
\end{abstract}

\textbf{Keywords:} finitely generated groups, classification.

\section{Introduction}


For a given group $X$, we shall denote by $Z(X)$ its center. Groups $G$ such that $G/Z(G) \simeq C_2 \times C_2$, where $C_2$ denotes a cyclic group of order 2, are called SLC groups and play a very important role in the construction of RA loops; i.e. loops such that loop ring over a commutative and associative ring with unity and no $2$-torsion is alternative but not associative. These loops were introduced by E.G. Goodaire \cite{good}, who proved that if $RL$ is alternative over one  ring $R$ as above, then it is also alternative over {\em all} such rings. Later, O. Chein and E.G. Goodaire \cite{CG} gave a means to construct all RA loops from SLC groups. For a  history of the subject, see also \cite{history}. More recently, SLC groups have played an important role in the study of involutions in group rings.

Using ideas from \cite{isomorfismo}, E. Jespers,
G. Leal and C. Polcino Milies   \cite{casofinito} gave a full description of finite
SLC groups  and classified all finite indecomposable RA loops,
up to isomorphism. These results are included in \cite[Chapter V]{livro}.

In this paper we consider the class of all finitely generated
groups $G$ such that $G/Z(G) \simeq C_p \times C_p$, where $C_p$
denotes a cyclic group of prime order $p$.
In section 2 we give a full description of such groups in terms of direct
factors and, in section 3, we classify indecomposable groups of this type in to nine families and show that groups in different families are not isomorphic.

As a consequence of this work, considering the case when $p=2$, it is possible to classify
finitely generated indecomposable RA loops. These results will be published elsewhere.

\section{A Decomposition Theorem}

Throughout this paper, $G$ will always denote a finitely generated group such that
$G/Z(G) \simeq C_p \times C_p$, where $C_p$
is a cyclic group of  order $p$.
In this section, we shall show that such a group decomposes in a very particular way.

\begin{theorem}\label{decomposition} A finitely generated group $G$ is such that $G / Z(G) \simeq C_p \times
C_p$ if and only if it can be written in the form $G = D \times A$
where $A$ is a finitely generated abelian group and $D$ is an  indecomposable group such that $D = \langle x,y,Z(D)\rangle$ where
$Z(D)$ is of the form
 $Z(D) = \langle t_1\rangle \times \langle z_2\rangle \times \langle z_3\rangle$, with:

 $(i)$  $o(t_1)= p^{m_1}, m_1 \geq 1$ and $s = [x,y]=x^{-1}y^{-1}xy = t_1^{p^{m_1-1}}$,

 $(ii)$ either $o(z_i)=
 p^{m_i}$ with $m_i \geq 0$ or $o(z_i) = \infty$, for $i=2,3$,

 $(iii)$
$x^p, y^p \in Z(D)$ with $x^p \in \langle t_1\rangle
\times \langle z_2\rangle$.
\end{theorem}

\emph{Proof:} Since $G / Z(G) \simeq C_p \times C_p$ there exist
elements $x,y \in G$ such that $G \simeq \langle x,y,Z(G) \rangle$ and $Z(G)$ is
 finitely generated.

 We can write $Z(G) = B
\times C \times F $ where $B$ and $C$ are finite, $B$ is a $p$-group, $p$ does not divide $|C|$ and $F$ is a direct product of a finite number of cyclic
groups of infinite order. As in [\cite{isomorfismo}, Theorem 1.2]
we can assume that $x^p,y^p \in B \times F$ and find a
decomposition  $B = \langle t_1\rangle \times \langle t_2\rangle \times \ldots \times
\langle t_k\rangle$, $o(t_i) = p^{m_i}$, such that $s=[x,y] = t_1^{p^{m_1-1}}$. 

Write $x^p = t_1^{a_1}t_2^{a_2} \ldots t_k^{a_k}l$ with $l \in F$.
If $p$ divides $a_i$, for some index $i$,
we consider
$x' = xt_i^{-a_i/p}$ so ${x'}^p$ has no component in
$\langle t_i\rangle$ and we still have that $G=\langle x', y, Z(G)\rangle$. Repeating this process, we can assume that $x^p=
t_1^{a_1}t_2^{a_2}\ldots t_m^{a_m}l$ where $\langle t_i^{a_i}\rangle=\langle t_i\rangle$ and
$o(t_2) \geq o(t_3) \geq \ldots \geq o(t_m)$

If $l=1$, then $x^p= t_1^{a_1}t_2^{a_2}\ldots t_m^{a_m}$ and we also have that
\[\langle t_2\rangle \times \langle t_3\rangle \times \ldots \times \langle t_m\rangle = \langle t_2^{a_2}\rangle \times \langle t_3^{a_3}\rangle \times \ldots \times \langle t_m^{a_m}\rangle
= \] \[ = \langle t_2^{a_2}t_3^{a_3}\ldots t_m^{a_m}\rangle \times \langle t_3^{a_3}\rangle
\times \ldots \times \langle t_m^{a_m}\rangle \]

Changing $t_2$ to
$t_2^{a_2}\ldots t_m^{a_m}$, we obtain that $x^p \in \langle t_1\rangle \times
\langle t_2\rangle$.

If $l \not=1$, then $x^p= t_1^{a_1}t_2^{a_2}\ldots t_m^{a_m}l$. We set
$l'= t_2^{a_2} \ldots t_m^{a_m}l$ and we get  that $x^p \in
\langle t_1\rangle \times \langle l'\rangle$, $\langle l'\rangle \subseteq F$.

There is a basis $\{ u_1,
u_2, \ldots , u_r\}$ of $F$ and a basis $\{v\}$ of
$\langle l'\rangle$ such that $v = u_1^{\alpha}$ for some $\alpha \in \mathbb{N}$.

We can write
\[ Z(G) = \langle t_1\rangle \times \langle t_2\rangle \times \ldots \times \langle t_k\rangle \times
\langle u_1\rangle \times \langle u_2\rangle \times \ldots \times \langle u_r\rangle \times C\]
and  also
\[ Z(G) = \langle t_1\rangle \times \langle u_1\rangle \times \langle t_2\rangle \times \ldots \times
\langle t_k\rangle \times \langle u_2\rangle \times \ldots \times \langle u_r\rangle \times C \]
obtaining that
$x^p \in \langle t_1\rangle \times \langle u_1\rangle$, with $o(u_1)= \infty$.

In any case, $x^p$ belongs, at most, to the product of the first two
factors, $\langle t_1\rangle \times \langle z_2\rangle$.\\

Write $y^p = t_1^{b_1}z_2^{b_2}t_3^{b_3} \ldots
t_k^{b_k}h$, with $h\in F$.

 As above, we can assume that $y^p =
t_1^{b_1}z_2^{b_2}t_3^{b_3} \ldots t_n^{b_n}h$ where
$\langle t_i^{b_i}\rangle = \langle t_i\rangle$, $3 \leq i\leq n$ . Reordering, if
necessary, we can also assume that $o(t_3) \geq o(t_4) \geq \ldots
\geq o(t_n)$.

If $h=1$ then $y^p= t_1^{b_1}z_2^{b_2}t_3^{b_3}\ldots t_n^{b_n}$
and
\[ \langle t_3\rangle \times \langle t_4\rangle \times \ldots \times \langle t_n\rangle = \langle t_3^{b_3}\rangle
\times \langle t_4^{b_4}\rangle \times \ldots \times \langle t_n^{b_n}\rangle = \]
\[ = \langle t_3^{b_3}t_4^{b_4}\ldots t_n^{b_n}\rangle \times \langle t_4^{b_4}\rangle \times
\ldots \times \langle t_n^{b_n}\rangle \]

Changing $t_3$ to
$t_3^{b_3}t_4^{b_4}\ldots t_n^{b_n}$ we can conclude that $y^p \in
\langle t_1\rangle \times \langle z_2\rangle \times \langle t_3\rangle$.

If $h \not=1$, then we set $ h' = t_3^{b_3}t_4^{b_4}\ldots t_n^{b_n}h$, where $y^p \in \langle t_1\rangle \times \langle z_2\rangle \times \langle h'\rangle$.

There exists a basis $\{u_1, u_2, \ldots, u_r\}$ of $F$ and a
basis $\{v\}$ of $\langle h'\rangle$ such that $v = u_1^{\beta}$ for some
$\beta \in \mathbb{N}$. Then, $y^p \in \langle t_1\rangle \times \langle z_2\rangle \times \langle u_1\rangle$ and,
reordering,
we see that $y^p$ belongs at most to the
first three factors of $Z(G)$.

In any case, we see that $y^p \in \langle t_1\rangle \times \langle z_2\rangle \times \langle z_3\rangle$ where $t_1$ is a $p$-element and $z_2, z_3$ are either $p$-elments or of infinite order.

Write $Z(G) = \langle t_1\rangle \times \langle z_2\rangle \times \langle z_3\rangle \times A$. Then

\[ G = \langle x,y,\langle t_1\rangle \times \langle z_2\rangle \times \langle z_3\rangle\rangle \times A = D \times A\] where $A$ is a finitely generated
abelian group.

The converse is straightforward. $\hfill{\fimprova}$
\\

We will need the following technical lemma, which is very similar to \cite[Lemma 3.1]{casofinito} and will be frequently used implicitly in  what follows.

\begin{lemma}\label{lema} With the notations above, the elements $x,y\in D$ can be chosen in such a way that  $x^p = t_1^{\alpha}z_2^{\beta}$ and
$y^p = t_1^{\gamma}z_2^{\epsilon}z_3^{\delta}$ with $\alpha, \beta, \gamma, \delta, \epsilon \in
\{0,1,\ldots, p-1\}$. \end{lemma}

\textbf{Proof.} Let $D = \langle x,y, Z(D)\rangle$ and set

\[x' = \left\{%
\begin{array}{ll}
xt_1^{\frac{p^{m_1}-\alpha}{p}}, & \hbox{if $\alpha \equiv 0$ mod $p$;} \\  \\

xt_1^{p^{m_1}-q}, & \hbox{if $\alpha = pq + r$, $r \in \{1,\ldots,p-1\}$.} \\
\end{array}%
\right. \]

Accordingly, we obtain ${x'}^p = z_2^{\beta}$ or ${x'}^p
= t_1^rz_2^{\beta}$. If $o(z_2) = p^{m_2}$, similar computations
will allow us to handle the exponent $\beta$ to obtain a new
element $x'$ such that ${x'}^p$ can be written as stated. If
$o(z_2) = \infty$, consider
\[x' = \left\{%
\begin{array}{ll}
xz_2^{\frac{-\beta}{p}}, & \hbox{if $\beta \equiv 0$ mod $p$;} \\ \\
xz_2^{-q}, & \hbox{if $\beta = pq+r$, $r \in \{1,\ldots, p-1\}$.} \\
\end{array}%
\right.\]

 Repeating this process for the exponents of $y^p$ we can
choose new elements $x,y \in G$ such that $\alpha, \beta, \gamma,
\delta, \epsilon \in \{0,1,\ldots, p-1\}$. $\hfill{\fimprova}$

\section{The Classification}

In view of the Theorem \ref{decomposition}, if $G$ is a finitely generated indecomposable group such that $G/Z(G) \simeq C_p \times C_p$ then it is of the form $G = \langle x,y,Z(G) \rangle$ with $x^p,y^p \in Z(G)$ with $rank[Z(G)] \leq 3$.

\begin{theorem}\label{posto1}
If $G$ is indecomplosable and $rank[Z(G)]=1$, then $G$ belongs to one of the following two families of groups:
\[ \mathcal{G}_1: \mbox{ groups with presentation } \langle x,y,t_1\,\, |\,\, x^p=y^p=t_1^{p^{m_1}}=1 \rangle.\]
\[ \mathcal{G}_2: \mbox{ groups with presentation }  \langle x,y,t_1 \,\,|\,\, x^p=y^p=t_1, t_1^{p^{m_1}}=1 \rangle. \]
where, in each case, $t_1$ is central, $m_1 \geq 1$ and
$s=[x,y]=t_1^{p^{m_1-1}}$.
\end{theorem}

\textbf{Proof.} Write $Z(G)=\langle t_1\rangle$. According to Theorem~\ref{decomposition} we have that $s =
t_1^{p^{m_1-1}}$ and $s^p =1$. In this case, $G$ is finite and we can write $x^p = t_1^a$, $y^p = t_1^b$ with $a,b \in \{0,1,\ldots, p-1\}$.

If $p =2$, then our statement coincides with \cite[Proposition 3.2]{casofinito}.

Suppose $p \neq 2$. If $a=b=0$ then $G \in \mathcal{G}_1$ and if $a=b=1$, $G \in \mathcal{G}_2$. If $a=0$ and $b \not=0$, change $x$ to $x' = xy$ so that $(x')^p = t_1^b$.  If $b =1$ we have that $G \in \mathcal{G}_2$ and if $b \in \{2,3, \ldots, p-1\}$, first we change $x$ to $x'= xy^n$ where $(n+1)b \equiv 1$ (mod $p$) and we can assume, by Lemma~\ref{lema}, that $(x')^p= t_1$.

If  we now change $y$ to $y'=(x')^my$ where $m + b \equiv 1$ (mod $p$), we see that $G = \langle x',y',Z(G) \rangle$ with $(x')^p=(y')^p = t_1$, which implies that $G \in \mathcal{G}_2$. $\hfill{\fimprova}$

\bigskip

\begin{proposition} \label{indcaso1}
Groups in any of the two families of Theorem \ref{posto1} are  indecomposable and groups in different families  are non-isomorphic.
\end{proposition}

\textbf{Proof.} Let $G \in \mathcal{G}_1$ or $\mathcal{G}_2$. It follows from Theorem~\ref{decomposition}, that $G = D \times A$ with $D$  indecomposable  and $A$ an abelian group. Hence, $Z(G) = Z(D) \times A$, and since $rank[Z(G)]=1$, we have that $A = \{1\}$ so $G$ is indecomposable group.\\

 Let $G_1 \in \mathcal{G}_1$ and $G_2 \in \mathcal{G}_2$ with $|G_1| = |G_2|.$ It is easily seen that $exp[G_1]=p^{m_1}$ and $exp[G_2]=p^{m_1+1}$, which shows that these groups are not isomorphic. $\hfill{\fimprova}$

\bigskip

\begin{theorem}\label{posto2} If $G$ is indecomposable and $rank[Z(G)]=2$, then $G$ belongs to one of the following four families of groups:
\begin{eqnarray*}
\mathcal{G}_3: & \mbox{ groups with presentation } & \langle x,y,t_1,t_2 \,\,|\,\, t_1^{p^{m_1}}=t_2^{p^{m_2}}=1, x^p=1,
y^p=t_2\rangle. \\
\mathcal{G}_4: & \mbox{ groups with presentation } & \langle x,y,t_1,t_2 \,\,|\,\, t_1^{p^{m_1}}=t_2^{p^{m_2}}=1, x^p=t_1, y^p=t_2 \rangle.\\
 \mathcal{G}_5:  & \mbox{ groups with presentation } & \langle x,y,t_1,u_1 \,\,|\,\,t_1^{p^{m_1}}=1, x^p = 1,
y^p = u_1 \rangle.\\
 \mathcal{G}_6: & \mbox{ groups with presentation } & \langle x,y,t_1,u_1 \,\,|\,\, t_1^{p^{m_1}}=1, x^p = t_1, y^p=u_1 \rangle.
\end{eqnarray*}
where, in each case, $t_1,t_2,u_1$ are central, $m_1, m_2 \geq 1$
and $s=[x,y] = t_1^{p^{m_1-1}}$.
\end{theorem}

\textbf{Proof.} Let $G$ be indecomposable. If $rank[Z(G)]=2$, then we can write $G = \langle x,y,Z(G) \rangle $ where $Z(G)= \langle t_1\rangle \times \langle z_2\rangle$ with $o(t_1)= p^{m_1}$ and $o(z_2)=p^{m_2}$ or $\infty$.

Let $\bar{G} = G/\langle z_2 \rangle = \langle \bar{x}, \bar{y}, \bar{t_1} \rangle$. Then $\bar{x}^p, \bar{y}^p \in \langle \bar{t_1} \rangle = Z(\bar{G})$ and, using Theorem~\ref{decomposition}, it is easily seen that $\bar{G}$ is indecomposable. According to Theorem \ref{posto1}, we have two possibilities: either
(1) $\bar{x}^p = \bar{y}^p= \bar{1}$  or (2) $\bar{x}^p = \bar{y}^p = \bar{t_1}$.

In the first case, we have four possibilities for $x^p$ and $y^p$:
\begin{itemize}
\item[(1.1)] $x^p = y^p = 1$
\item[(1.2)] $x^p = 1$, $y^p = z_2^a$
\item[(1.3)] $x^p=z_2^a$, $y^p =1$
\item[(1.4)] $x^p = z_2^a$, $y^p=z_2^b$
\end{itemize}

Recall that, because of Lemma \ref{lema}, we can assume that $a,b \in \{1,2, \ldots, p-1\}$.

In $(1.1)$, we get $G = \langle x,y,t_1 \rangle \times \langle z_2 \rangle$, which is not indecomposable.

Assume first that $z_2 = t_2$ with $o(t_2) = p^{m_2}$ so $Z(G)$ is finite. We can also assume that $p \neq 2$, since the case $p=2$ is \cite[Proposition 3.3]{casofinito}.

In $(1.2)$, we can change $t_2$ to $t_2'=t_2^a$  and conclude that $G \in \mathcal{G}_3$. The case (1.3) is symmetric to the previous one and also defines a group in $\mathcal{G}_3$. In  case (1.4), change $x$ to $x'=x^ny$ where $n$ is such that $an + b \equiv 0$ (mod $p$). Using Lemma \ref{lema}, we can assume that $x'^p=1$ and $y^p=t_2^b$ and changing $t_2$ to $t_2^b$ we obtain that $G \in \mathcal{G}_3$.

Suppose now that $z_2 = u_1$ with $o(u_1)= \infty$. In  case $(1.2)$, if $a=1$, then $G \in \mathcal{G}_5$. If $p=2$, then by Lemma~\ref{lema}, we can assume that $a=1$  so $G \in \mathcal{G}_5$. If $p \not=2$, changing $y$ to $y'=y^b$ where $b$ is such that $ab \equiv 1$ (mod $p$), we can assume that $x^p=1$ and $y'^p=u_1$ so we get that $G \in \mathcal{G}_5$.

The  case (1.3) is symmetric and in  case (1.4), if $p=2$, then $x^2=y^2=u_1$. Changing $y$ to $y'=xy$, we can suppose that $y'^2=x^2y^2s=t_1^{2^{m_1-1}}$. If $m_1>1$, we have $G \in \mathcal{G}_5$ and if $m_1=1$, then $G \in \mathcal{G}_6$. If $p\neq 2$, changing $x$ to $x'=xy^n$ where $n$ is such that $a+nb \equiv 0$ (mod $p$), we obtain $G = \langle x' , y, Z(G) \rangle$ with $(x')^p=1$ and $y^p=u_1^b$ which implies, by  case $(1.2)$, that $G \in \mathcal{G}_5$.\\

If $\bar{x}^p = \bar{y}^p=\bar{t_1}$, we  have again four possibilities for $x^p$ and $y^p$:

\begin{itemize}
\item[(2.1)] $x^p = y^p = t_1$
\item[(2.2)] $x^p = t_1$, $y^p = t_1z_2^a$
\item[(2.3)] $x^p = t_1z_2^a$, $y^p = t_1$
\item[(2.4)] $x^p = t_1z_2^a$, $y_p = t_1z_2^b$
\end{itemize}

The first one can not occur, as else $G = \langle x,y,t_1 \rangle \times \langle z_2 \rangle$ is decomposable.

Once again, we shall assume first that $z_2 = t_2$ with $o(t_2)=p^{m_2}$. We can suppose that $p \neq 2$, because otherwise the result is \cite[Proposition 3.3]{casofinito}. In $(2.2)$, changing $y$ to $y'=x^ny$ with $n+1 \equiv 0$ (mod $p$), we can assume that $x^p=t_1$ and $(y')^p=t_2^a$ which implies $G \in \mathcal{G}_4$ if we change $t_2$ to $t_2^a$.

The case $(2.3)$ is symmetric and in  case (2.4), changing $x$ to $x'=x^{(p-1)}y$, as we can assume by Lemma~\ref{lema} that $t_1^p=1$, we get $x'^p = t_2^{a(p-1)+b}$.

Set $a'= a(p-1)+b$. If $a'\not\equiv 0$ (mod $p$), then we can change $y$ to $y'=x^my$ with $b+ma' \equiv 0$ (mod $p$) and obtain $x'^p=t_2^{a'}$, $y'^p = t_1$; and taking $t_2' = t_2^{a'}$ we see that $G \in \mathcal{G}_4$.

If $a' \equiv 0$ (mod $p$), then we can suppose that $x^p=1$, $y^p=t_1t_2^b$. If $m_2 \geq m_1$, then $Z(G) = \langle t_1 \rangle \times \langle t_2 \rangle = \langle t_1 \rangle \times \langle t_2^b \rangle = \langle t_1 \rangle \times \langle t_1t_2^b \rangle$. Changing $t_2$ to $t_2'=t_1t_2^b$, we obtain $G \in \mathcal{G}_3$. If $m_2 < m_1$, we set $t_1'= t_1t_2^b$ and then $\langle t_1 \rangle \times \langle t_2 \rangle = \langle t_1t_2^b \rangle \times \langle t_2 \rangle$, where $s = (t_1t_2^b)^{p^m_1-1}$ which implies that $G = \langle x,y,t_1' \rangle \times \langle t_2 \rangle$ is decomposable, a contradiction.

Finally, suppose that $z_2 = u_1$ with $o(u_1) = \infty$. In  case $(2.2)$, changing $y$
to $y'=y^b$ with $ab \equiv 1$ (mod $p$), we can assume that $x^p = t_1$ and
$y'^p=t_1^bu_1$. Changing $u_1$ to $u_1'=t_1^bu_1$ we have that $G \in \mathcal{G}_6$ and
the case $(2.3)$ is symmetric to this one.

 In  case $(2.4)$, if $p=2$, we have $x^2=y^2=t_1u_1$ so
changing $x$ to $x'=xy$, we see that $x'^2=t_1^{2^{m_1-1}+2}$.
If $m_1>1$ we can suppose that $x'^2 = 1$ and $y^2 = t_1u_1$ and setting $u_1'=t_1u_1$
we see that $G \in \mathcal{G}_5$. If $m_1 =1$, we have $x'^2 = t_1$,
$y^2 = t_1u_1$ so, if we change $u_1$ to $u_1'=t_1u_1$, we see that $G \in \mathcal{G}_6$ .

Now, if $p\neq 2$ and $x^p = t_1u_1^a$, $y^p =t_1u_1^b$, we change $x$ to $x'=xy^{p-1}$. Then $x'^p=u_1^{a+(p-1)b}$.
Write $a'=a+(p-1)b$. If $a' \not\equiv 0$ (mod $p$), we change $y$ to $y'=(x')^my$ where $b+ma' \equiv 0$ (mod $p$) so $y'^p = t_1$ and $x'^p = u_1^{a'}$. Then, change $x'$ to $(x')^c$ where $a'c \equiv 1$ (mod $p$). We have obtained  $x', y' \in G$ such that $x'^p = u_1$, $y'^p = t_1$ which implies that $G \in \mathcal{G}_6$. If $a' \equiv 0$ (mod $p$), there exists $x' \in G$ with $x'^p=1$. Changing $y$ to $y'=y^d$ with $db \equiv 1$ (mod $p$), we obtain  that $G \in \mathcal{G}_5$ if we also change $t_1^du_1$ by $u_1'$. $\hfill{\fimprova}$

\bigskip
In order to study the non decomposability of groups in these families, we will need the following result.

\begin{lemma}\label{lema2}
Let $G = \langle x,y,Z(G) \rangle = D \times A$ where $x^p,y^p \in Z(G)$ and $D$ and  $A$ are as in Theorem~\ref{decomposition}. Then, there exist $x_1, y_1 \in D$  and $z,z' \in Z(G)$ such that $x= x_1z^{-1}$ and $y = y_1z^{-1}$.
\end{lemma}

\textbf{Proof.}\, Since $G / Z(G) \simeq C_p \times C_p$ and $|C_p \times C_p| = p^2$, we can write $G$ as a disjoint union of cosets of the form $x^iy^jZ(G)$ with $0 \leq i,j \leq p-1$. We have that $D$ contains elements of two different cosets since $D$ is non abelian. Let $x^ay^bz_1$ and $x^ny^mz_2 \in D$ with $z_1, z_2 \in Z(G)$. Then,
\[ (x^ay^bz_1)(x^ny^mz_2) = x^{a+n}y^{b+m}z_1z_2s^{bn}\]
\[ (x^ny^mz_1)(x^ay^bz_2) = x^{n+a}y^{m+b}z_1z_2s^{am}\]

Since these two elements do not commute, $am - bn \not= 0$ (mod $p$), which implies that the linear system below has a solution:
\[\left\{
    \begin{array}{ll}
      a\alpha + n\beta \equiv 0 \,\, (mod \,\, $p$) \\
      b\alpha + m\beta \equiv 1 \,\, (mod \,\, $p$)
    \end{array}
  \right.
\]
Hence, there exist $\alpha, \beta, \alpha', \beta' \in \mathbb{N}$ such that
\[y_1 = (x^ay^bz_1)^{\alpha}(x^ny^mz_2)^{\beta} = yz \,\,\, \hbox{and} \,\,\, y_2 = (x^ay^bz_1)^{\alpha '}(x^ny^mz_2)^{\beta'} = xz'\]
where $z, z' \in Z(G)$.
$\hfill{\fimprova}$

\begin{proposition} \label{indcaso2}
Groups in any of the four families of Theorem \ref{posto2} are  indecomposable and groups in different families  are non-isomorphic.
\end{proposition}

\textbf{Proof.} Let $G$ be a group in one of the four families of the previous theorem. By Theorem \ref{decomposition} we can write $G = D \times A$, where $A$ is an
 abelian group and $D$ is indecomposable. Then, $Z(G)= \langle t_1\rangle \times \langle z_2\rangle = Z(D)\times A$. If $A$ is non-trivial, again by Theorem~\ref{decomposition}, as $s \in Z(D)$ and $s^p=1$ we can conclude that $Z(D)$ is cyclic of order $p^r$.
 Suppose that
$Z(D)=\langle t_1^{a}z_2^{b}\rangle$ with
$o(t_1^{a}z_2^{b}) = p^r$. Then $s = t_1^{ap^{r-1}j}z_2^{bp^{r-1}j} = t_1^{p^{m_1}-1}$ where $p$ does not divide $j$.

If $G \in \mathcal{G}_3$ or  $G \in \mathcal{G}_4$ then $z_2 = t_2$ with $o(t_2)=p^{m_2}$. Hence, $p^{m_2} | bp^{r-1}j$ which implies that either $r-1 \geq m_2$ or $b$ is a multiple of $p$. By lemma \ref{lema2}, we know that there exists $y_1 \in D$ such that $y_1 = yz$ with $z \in Z(G)$. So, ${y_1}^p =
y^pz^p= t_2(t_1^{\alpha}t_2^{\beta})^p \in D \cap Z(G) \subseteq Z(D)$ and then we can write $y_1^p=(t_1^at_2^b)^n$ so, $t_1^{na}=t_1^{\alpha p}$ and
 $t_2^{\beta p +1}=t_2^{nb}$.
As a consequence, we see that $p$  divide neither $n$ nor $b$ so, from the above, it follows that $r-1\geq m_2$.

As $Z(D)$ is a direct factor of $Z(G)$ of order $p^r$ and the factors of $Z(G)$ have orders $p^{m_1}$ and $p^{m_2}$ we can conclude that $r = m_1$.

If $a$ is a multiple of $p$, we have that $(t_1^at_2^b)^{p^{m_1}-1}=1$,
which is a contradiction, as $o(t_1^at_2^b) = p^r = p^{m_1}$. Also, since
$p^{m_1} | (\alpha p - an)$, we have that $p$ divide $an$ which implies that $p | a$ (since $p \nmid n$), again a contradiction.

If $G \in \mathcal{G}_5$ or $\mathcal{G}_6$, then $z_2 = u_1$ with $o(u_1)= \infty$.
It follows that $u_1^{bp^{r-1}j}=1$ and hence $b =0$, $Z(D)= \langle t_1^a \rangle$ and $t_1^{an} = t_1^{\alpha p}u_1^{1+\beta p}$. So, we would have $1 + \beta p = 0$ which is a contradiction since $\beta \in \mathbb{Z}$.

\bigskip
 Let $G_i \in \mathcal{G}_i$ for $i = 3,4,5,6$. Because of cardinality considerations, we have only to prove that $G_3 \not\cong G_4$ and $G_5 \not\cong G_6$. To prove that $G_3 \not\cong G_4$, note that $x$ is a non central element of order $p$ in $G_3$. We will show that there does not exist such an element in $G_4$. Suppose that $w = x^ay^bt_1^ct_2^d \in G_4$, $w \not\in Z(G_4)$ and $w^p=1$. Because of \cite[Proposition 3.3]{casofinito},  we can assuma $p\neq 2$. Hence $w^p = t_1^{a+cp}t_2^{b+dp} = 1$ and as $\langle t_1 \rangle \cap \langle t_2 \rangle = \{1\}$, we have that $p | a$ and $p | b$ which implies that $w \in Z(G_4)$, a contradiction. A similar argument shows that $G_5 \not\cong G_6$.$\hfill{\fimprova}$

\bigskip

Finally, the next theorem shows what happens in the case of the groups of rank  three.

\begin{theorem}\label{posto3}
If $G$ is indecomposable and $rank[Z(G)]=3$, then $G$ belongs to one of the following three families of groups:
\begin{eqnarray*}
 \mathcal{G}_7: & \mbox{ groups with presentation } &  \langle x,y,t_1,t_2,t_3 \,\,|\,\,
t_1^{p^{m_1}}=t_2^{p^{m_2}}=t_3^{p^{m_3}}=1, x^p=t_2, y^p=t_3\rangle.\\
 \mathcal{G}_8: & \mbox{ groups with presentation } & \langle x,y,t_1,t_2,u_1 \,\,|\,\, t_1^{p^{m_1}} = t_2^{p^{m_2}} =1, x^p=t_2, y^p=u_1
\rangle. \\
 \mathcal{G}_9: & \mbox{ groups with presentation } &   \langle x,y,t_1,u_1,u_2 \,\,|\,\, t_1^{p^{m_1}}=1, x^p=u_1, y^p= u_2\rangle.
 \end{eqnarray*}
where, in each case, $t_1, t_2,t_3,u_1,u_2$ are central, $m_1,
m_2,m_3 \geq 1$ and $s= t_1^{p^{m_1-1}}$.
\end{theorem}

\textbf{Proof.} If $G$ is indecomposable and $rank[Z(G)]=3$, by Theorem \ref{decomposition} we have that
$Z(G)=\langle t_1 \rangle \times \langle z_2 \rangle  \times  \langle z_3 \rangle$, $o(t_1) = p^{m_1}$, $m_1 \geq 0$,
and $o(z_i) = p^{m_i}$ or $\infty$, $i = 2,3$.

Let $\overline{G} = G / \langle z_3 \rangle = \langle \bar{x}, \bar{y}, \bar{t_1}, \bar{z_2} \rangle$ where $\bar{x}^p, \bar{y}^p \in \langle \bar{t_1} \rangle \times \langle \bar{z_2} \rangle = Z(\overline{G})$. If $\overline{G}$ is indecomposable, by Theorem \ref{posto2}, we have two possibilities: either (1) $\bar{x}^p = \bar{1}, \bar{y}^p = \bar{z_2}$ or (2) $\bar{x}^p = \bar{t_1}^p, \bar{y}^p = \bar{z_2}$.

If (1) occurs, we have four possibilities for $x^p, y^p$:
\begin{itemize}
\item[(1.1)] $x^p=1, y^p = z_2$
\item[(1.2)] $x^p=1, y^p = z_2z_3^a$
\item[(1.3)] $x^p= z_3^a, y^p z_2$
\item[(1.4)] $x^p=z_3^a, y^p = z_2z_3^b$
\end{itemize}
where $a,b \in \{1,2,\ldots, p-1\}$.

In  case (1.1), $G = \langle x,y,t_1,z_2 \rangle \times  \langle z_3 \rangle$ is always decomposable.

Suppose $z_i = t_i$ with $o(t_i) = p^{m_i}$, $i=2,3$.

In  case $(1.2)$, if $m_2 \geq m_3$, then $Z(G) = \langle t_1 \rangle \times \langle t_2t_3^a \rangle \times \langle t_3 \rangle$ and setting $t_2'= t_2t_3^a$, we have that $G$ is decomposable, a contradiction. If $m_2 < m_3$, changing $t_3$ to $t_3' = t_2t_3^a$ we can conclude that $G = \langle x,y t_1, t_3' \rangle \times \langle t_2 \rangle$, a contradiction. So, under these conditions,  case $(1.2)$ can not occur.

In  case $(1.3)$, setting $t_3^a = t_3'$, we have that $G \in \mathcal{G}_7$.

In the last case, if $p =2$, \cite[Proposition 3.4]{casofinito} shows that $G \in \mathcal{G}_7$. If $p\not= 2$, changing $y$ to $y'= x^ny$ with $b + an \equiv 0$ (mod $p$), we have that $(y')^p = t_2$. Setting $t_3^a = t_3'$, we conclude that $G \in \mathcal{G}_7$.

Now, suppose $z_2 = t_2$, with $o(t_2)= p^{m_2}$ and $z_3 = u_1$ with $o(u_1) = \infty$. In $(1.2)$, if $p=2$, $x^2=1$ and $y^2= t_2u_1$, changing $u_1$ to $u_1' = t_2u_1$ we have $G = \langle x,y,t_1,u_1' \rangle \times \langle t_2 \rangle$, which is a contradiction. If $p \neq 2$, changing $y$ to $y'= y^b$, with $ab \equiv 1$ (mod $p$) and changing $u_1$ to $u_1' = t_2^bu_1$ we have that   $(y')^p = t_2^bu_1$, which implies that $G$ is  decomposable.

The case $(1.3)$, if we set $x'=x^b$ with $ab \equiv 1$ (mod $p$), we have that $G \in \mathcal{G}_8$.

In the last case, if $p=2$, we have $x^2 = u_1$ and $y^2 = t_2u_1$ so changing $y$ to $y'= xy$, we may suppose that $y'^2= t_1^{2^{m_1-1}}t_2$. If $m_1 =1$, we change $t_2$ to $t_2' = t_1t_2$, and see that $G \in \mathcal{G}_8$. If $m_1 > 1$ we can suppose $y'^2 = t_2$, and hence $G \in \mathcal{G}_8$.

Now, if $p\neq 2$, setting $y' = x^ny$ where $an + b \equiv 0$ (mod $p$) and $x' = x^c$ where $ac \equiv 1$ (mod $p$) we can assume that $y'^p = t_2$ and $x'^p = u_1$, so $G \in \mathcal{G}_8$.

Finally, suppose that $z_i = u_{i-1}$,  with $o(u_{i-1}) = \infty$, $i = 2,3$. In  case $(1.2)$, changing $u_1$ to $u_1' = u_1u_2^a$, we obtain $G = \langle x,y, t_1, u_1' \rangle \times \langle u_2 \rangle$ which is decomposable. In case $(1.3)$, changing $x$ to $x' = x^b$ with $ab \equiv 1$ (mod $p$), we have $G \in \mathcal{G}_9$. In the last case, working as in $(1.3)$ and renaming $u_1' = u_1u_2^a$, we have $G = \langle x', y, \langle t_1 \rangle, \times \langle u_1u_2^a \rangle \times \langle u_2 \rangle \rangle$ with $x'^p= u_2$ and $y^p = u_1$ which implies that $G \in \mathcal{G}_9$.\\

If $(2)$ occurs, then $\bar{x}^p = \bar{t_1}$, $\bar{y}^p=\bar{z_2}$ and we have again four possibilities for $x^p, y^p$:
\begin{itemize}
\item[(2.1)] $x^p = t_1, y^p = z_2$
\item[(2.2)] $x^p = t_1, y^p = z_2z_3^a$
\item[(2.3)] $x^p = t_1z_3^a, y^p = z_2$
\item[(2.4)] $x^p = t_1z_3^a, y^p = z_2z_3^b$
\end{itemize}
where $a,b \in \{1,2, \ldots, p-1\}$.

In the first case, $G = \langle x,y,t_1,z_2 \rangle \times \langle z_3 \rangle$ is decomposable, a contradiction.

Suppose that $z_i = t_i$ with $o(t_i) = p^{m_i}, i = 2,3$. In case $(2.2)$, if $m_2 \geq m_3$, we change $t_2$ to $t_2' = t_2t_3^a$ and otherwise, we change $t_3$ by $t_3' = t_2t_3^a$. In both cases, $G$ is decomposable.

In case $(2.3)$, if $m_3 \geq m_1$, we can change $t_3$ to $t_3' = t_1t_3^a$ and we obtain that $G \in \mathcal{G}_7$. If $m_3 < m_1$, we take $t_1' = t_1t_3^a$ (this is possible since $s = t_1^{p_{m_1-1}} = (t_1t_3^a)^{p^{m_1-1}}$ and obtain $G = \langle x,y, t_1', t_2 \rangle \times \langle t_3 \rangle$, a contradiction.

In the last case, we can suppose $p \not=2$ since in the case $p =2$ it was shown in \cite[Proposition 3.3]{casofinito} that $G \in \mathcal{G}_7$. If $m_1 \geq m_2 \geq m_3$, then changing $x$ to $x'= xy^n$ where $a + nb \equiv 0$ (mod $p$) we have $x'^p= t_1t_2^n$. If $n \equiv 0$ (mod $p$), we can suppose that $x'^p = t_1$. In this case, $G$ is decomposable, since we can change $t_2$ to $t_2' = t_2t_3^b$ and obtain $G = \langle x,y,t_1,t_2' \rangle \times \langle t_3 \rangle$. So, assume that $x'^p = t_1t_2^a, y^p = t_2t_3^b$ where $a,b \in \{1,2, \ldots, p-1\}$. We can write $Z(G) = \langle t_1t_2^a \rangle \times \langle t_2t_3^b\rangle \times \langle t_3 \rangle$ and see that $G$ is decomposable changing $t_1$ to $t_1' = t_1t_2^a$ and $t_2$ to $t_2' = t_2t_3^b$. If $m_1 \geq m_2$ and $m_2 < m_3$, we change again $x$ to $x' = xy^n$ where $a+nb \equiv 0$ (mod $p$) and if $n\not\equiv 0$ (mod $p$), we take $t_1' = t_1t_2^a$ and $t_3' = t_2t_3^b$. In this way, $Z(G) = \langle t_1' \rangle \times \langle t_2 \rangle \times \langle t_3' \rangle$, hence $G = \langle x', y,t_1', t_3' \rangle \times \langle t_2 \rangle$, a contradiction. If $n \equiv 0$ (mod $p$), we have $x'^p = t_1, y^p = t_2t_3^b$, which implies $G$ decomposable changing $t_3$ to $t_3' = t_2t_3^b$. If $m_3 < m_1 < m_2$, since $x^p = t_1t_3^a, y^p = t_2t_3^b$, we set $t_1' = t_1t_3^a$ and $t_2' = t_2t_3^b$ and obtain that $G$ is decomposable. Suppose that $m_1 \geq min\{m_2,m_3\}$. Changing $x$ to $x' = xy^n$ where $a+nb \equiv 0$ (mod $p$), we obtain $x'^p = t_1t_2^n$. If $n \equiv 0$ (mod $p$) we have seen that this implies $G$ decomposable; so we can assume that $n \not\equiv 0$ (mod $p$) and change $y$ to $(x')^my$ with $mn + 1 \equiv 0$ (mod $p$). Then, we have $x'^p = t_1t_2^n, y'^p = t_1^mt_3^b$ and, since $Z(G) = \langle t_1 \rangle \times \langle t_1t_2^n \rangle \times \langle t_1^mt_3^b \rangle$, we can conclude that $G \in \mathcal{G}_7$ changing $t_2$ to $t_2' = t_1t_2^n$ and $t_3$ to $t_3' = t_1^mt_3^b$.
\\

Now, let $z_2 = t_2$ with $o(t_2) = p^{m_2}$ and $z_3 = u_1$ with $o(u_1) = \infty$. In case $(2.2)$, changing $y$ to $y'=y^b$ where $ab \equiv 1$ (mod $p$) and setting $u_1'= t_2^bu_1$, we obtain $G = \langle x,y', t_1, u_1' \rangle \times \langle t_2 \rangle$ which is decomposable. In the third case, if we change $x$ to $x'=x^b$ with $ab \equiv 0$ (mod $p$) and set $u_1' = t_1^bu_1$ then we see that $G \in \mathcal{G}_8$.

 In case $(2.4)$, if $p=2$ we have $x^2= t_1u_1$ and $y^2 = t_2u_1$. Changing $x$ to $x' = xy$ we obtain $x'^2 =x^2y^2s = t_1^{2^{m_1-1}+1}t_2u_1^2$. If $m_1= 1$, we can assume that $x'^p= t_2$; hence, we obtain that $G \in \mathcal{G}_8$ setting $u_1'=t_2u_1$. If $m_1>1$, we can suppose that $x'^p = t_1t_2$. If $m_2 \geq m_1$ then $Z(G) = \langle t_1 \rangle \times \langle t_1t_2 \rangle \times \langle t_2u_1 \rangle$ and $G \in \mathcal{G}_8$. If $m_2 < m_1$, we can write $Z(G) = \langle t_1t_2 \rangle \times  \langle t_2 \rangle \times \langle t_2u_1 \rangle$, which implies that $G$ is decomposable. If $p \not= 2$, we change $x$ to $x' = xy^n$ with $a + nb \equiv 0$ (mod $p$) and $y$ to $y' = y^c$ with $bc \equiv 1$ (mod $p$). In this way, we can assume that $x'^p = t_1t_2^n$ and $y'^p = t_2^cu_1$. If $n \equiv 0$ (mod $p$) then $x'^p = t_1$ and setting $u_1' = t_2^cu_1$ it follows that $G = \langle x',y', t_1, u_1' \rangle \times \langle t_2 \rangle$. If $n \not\equiv 0$ (mod $p$) we recall that we can assume $ n \in \{1,2,\ldots, p-1 \}$. If $m_2 \geq m_1$ then $G = \langle x', y', \langle t_1 \rangle \times \langle t_1t_2^n \rangle \times \langle t_2^cu_1 \rangle \rangle \in \mathcal{G}_8$ and if $m_2 < m_1$ then $G = \langle x',y', \langle t_1t_2^n \rangle \times \langle t_2^cu_1 \rangle \rangle \times \langle t_2 \rangle$ is decomposable.
\\

If $z_i = u_{i-1}$ with $o(u_{i-1}) = \infty, i = 2,3$ then in case $(2.2)$ setting $u_1' = u_1u_2^a$, it follows that $G$ is decomposable. In case $(2.3)$, changing $x$ to $x' = x^b$ with $ab \equiv 1$ (mod $p$), we can assume that $x'^p = t_1^bu_2$ and $y^p = u_1$ which implies $G \in \mathcal{G}_9$ changing $u_2$ to $u_2' = t_1^bu_2$. In case $(2.4)$ changing $x$ to $x'= x^c$ where $c$ is such that $ac \equiv 1$ (mod $p$), we can assume that $x'^p = t_1^cu_2, y^p = u_1u_2^b$. Setting $u_1' = u_1u_2^b$ and $u_2' = t_1^cu_2$ we obtain that $G \in \mathcal{G}_9$.
\\

Now we turn to the case when $\overline{G}$ is decomposable. By Theorem \ref{decomposition}, $\overline{G} = \bar{D} \times \bar{A}$ where $\bar{A}$ is non trivial abelian group and $\bar{D} = \langle \bar{d_1}, \bar{d_2}, Z(\bar{D}) \rangle$ is indecomposable. Since $\langle z_3 \rangle \subseteq D \cap Z(G) \subseteq Z(D)$ we have that $ G = D.A$ where $D \cap A = \langle z_3 \rangle$ and $ D = \langle d_1, d_2, Z(D) \rangle$. Hence we can write $G = \langle d_1, d_2, Z(G) \rangle$. We claim that $Z(D)\subseteq Z(G)$. In fact, otherwise there would exist an element $d \in Z(D)$ such that $d \not\in Z(G)$ and then $s = d^{-1}a^{-1}da$ for some $a \in A$. Since $D$ and $A$ are normal subgroups of $G$ , we would have $s = t_1^{p^{m_1-1}} \in D \cap A = \langle z_3 \rangle$, a contradiction.

Let $H = [\langle t_1 \rangle \times \langle t_2 \rangle ] \cap Z(D)$. Note that $H \not= \emptyset$ as $ s = t_1^{p^{m_1-1}} \in Z(D)$. Since $Z(D) \subseteq Z(G) = \langle t_1 \rangle \times \langle z_2 \rangle \times \langle z_3 \rangle $ and $\langle z_3 \rangle \subset Z(D)$, we have that $Z(D) = H \times \langle z_3 \rangle$ and $H \simeq Z(D)/\langle z_3 \rangle \subseteq Z(\bar{D})$. Since $Z(\bar{D})$ is cyclic, we conclude that $Z(D)$ is of rank 2. Write $Z(D) = \langle t_1^az_2^b \rangle \times \langle z_3 \rangle$. We have two possibilities: either $d_1^p = 1, d_2^p = z_3$ or $d_1^p = t_1^az_2^b, d_2^p = z_3$. The first one implies that $G$ is decomposable. In the second one, using Lemma \ref{lema} we can suppose $a,b \in \{1,2, \ldots, p-1 \}$. Since $s \in Z(D)$ and $s^p =1$ then $z_2 = t_2$ with $o(t_2) = p^{m_2}$. If $m_2 \geq m_1$, we change $t_2$ to $t_2' = t_1^at_2^b$. If $z_3 = t_3$ with $o(t_3) = p^{m_3}$ we obtain $G \in \mathcal{G}_7$ and if $z_3 = u_1$ with $o(u_1) = \infty$ we obtain that $G \in \mathcal{G}_8$. If $m_2 < m_1$, note there exists a positive integer $\alpha $ such that $a\alpha \equiv p^{m_1-1}$ (mod $p^{m_1}$) and $b\alpha \equiv 0$ (mod $p^{m_2}$). Then $s = t_1^{p^{m_1-1}} = (t_1^at_2^b)^{\alpha}$ and we can change $t_1$ to $t_1' = t_1^at_2^b$ and obtain that $G$ is decomposable when $z_3 = t_3$ and also when $z_3 = u_1$. $\hfill{\fimprova}$

\begin{proposition} \label{indcaso3}
Groups in any of the three families of Theorem \ref{posto3} are indecomposable and groups in different families are non-isomorphic.
\end{proposition}

\textbf{Proof.} Suppose that a group $G = \langle x,y,Z(G)\rangle$ with $Z(G) = \langle t_1 \rangle \times \langle z_2 \rangle \times \langle z_3 \rangle$ in either $\mathcal{G}_7$, $\mathcal{G}_8$ or $\mathcal{G}_9$ is decomposable. By Theorem \ref{decomposition}, $G$ contains a unique non trivial commutator $s$ and we can write $G = D \times A$ where $A$ is an non trivial abelian group, $D = \langle d_1,d_2, Z(D) \rangle$ is indecomposable with $D/Z(D) \simeq C_p \times C_p$ and $s \in Z(D), s^p = 1$. As $Z(G) = Z(D) \times A$ and $rank[Z(G)]=3$, we have that $rank[Z(D)]=1$ or $2$.

 If $rank[Z(D)]=1$ write $Z(D) = \langle t_1^{a}z_2^bz_3^c \rangle$. Since $s \in Z(D)$ and $s^p=1$, we can conclude that $Z(D)$ is a finite cyclic group and, from Theorem \ref{posto1}, that $o(t_1^{a}z_2^bz_3^c) = p^r$ for some integer $r$. We know, by Lemma \ref{lema2}, that there exists $y_1 \in D$ such that $y_1 = yz$ with $z = t_1^{\alpha}z_2^{\beta}z_3^{\gamma} \in Z(G)$. So, $y_1^p \in D \cap Z(G) \subseteq Z(D)$ and we can write
\[ y_1^p = (t_1^{a}z_2^bz_3^c)^n = z_3t_1^{\alpha p}z_2^{\beta p}z_3^{\gamma p}.\]

If $G \in \mathcal{G}_7$, the equation above implies that $c$ and $n$ are not multiples of $p$. Lemma \ref{lema2} guarantees that there exists also $x_1 \in D$ such that $x_1 = xz'$ with $z' = t_1^{\alpha'}t_2^{\beta'}t_3^{\gamma'} \in Z(G)$. So, $x_1^p \in D \cap Z(G) \subseteq Z(D)$ and we can write

\[ x_1^p = (t_1^{a}z_2^bz_3^c)^{n'} = t_2t_1^{\alpha' p}t_2^{\beta' p}t_3^{\gamma' p}.\]
We conclude that also $b$ is not a multiple of $p$. Since $s \in Z(D)$ and $s^p=1$, we can write $s = (t_1^{a}t_2^bt_3^c)^{p^{r-1}j}$ where $p$ does not divide $j$. Hence,
\[ s = t_1^{p^{m_1-1}} = t_1^{ajp^{r-1}}t_2^{bjp^{r-1}}t_3^{cjp^{r-1}}\]
which implies that $p^{m_2} | (bjp^{r-1})$ and $p^{m_3} | (cjp^{r-1})$. Then $r-1 \geq m_2$ and $r-1 \geq m_3$ as $b,c$ and $j$ are not multiples of $p$. Since $Z(D)$ is a direct factor of $Z(G)$ of order $p^r$ and $Z(G)$ has three direct factors, each of order $p^{m_i}, i=1,2,3$, we have that $r = m_1$. If $a$ is a multiple of $p$, then $(t_1^{a}t_2^bt_3^c)^{p^{m_1-1}} =1$ which is a contradiction since $o(t_1^{a}t_2^bt_3^c)= p^r = p^{m_1}.$ But, we saw that $p^{m_1} | (an + \alpha p)$ which implies that $p| (an)$, again a contradiction since $a$ and $n$ are not multiples of $p$. Hence, if $G \in \mathcal{G}_7$, it is not decomposable.

If $G \in \mathcal{G}_8$ then $Z(D) = \langle t_1^at_2^b \rangle$ and $(t_1^at_2^b)^n = u_1t_1^{\alpha p}t_2^{\beta p}u_1^{\gamma p}$. Hence $u_1^{1+\gamma p} = 1$, a contradiction since $o(u_1) = \infty$. Similar arguments shows that if $G \in \mathcal{G}_9$ and $rank[Z(D]=1$, then $G$ is not decomposable.


Now assume that $Z(D) = \langle \alpha \rangle \times \langle \beta \rangle$. By  Theorem \ref{posto2}, we have
\begin{enumerate}
\item $d_1^p = 1$ and $d_2^p = \beta$ or
\item $d_1^p = \alpha$ and $d_2^p = \beta$
\end{enumerate}
Since $d_1 \in G = \langle x,y,t_1,z_2,z_3 \rangle$, we can write $d_1 = x^ry^st_1^mz_2^nz_3^q$ with $r,s,m,n,q \geq 0$. Then,
\begin{equation}\label{horrible}
d_1^p = (x^p)^r(y^p)^st_1^{mp}z_2^{np}z_3^{qp}(s^{\frac{p(p-1)}{2}})^{rs} = t_1^{mp}z_2^{r+np}z_3^{s+qp}(s^{\frac{p(p-1)}{2}})^{rs}
\end{equation}

If the first possibility occurs, we have $z_2^{r+np} = 1$ and $z_3^{s+qp} = 1$. If $G \in \mathcal{G}_7$ then $p^{m_2} | (r+np)$ and $p^{m_3} | (s+qp)$ so we can conclude that $r$ and $s$ are multiples of $p$ and hence $d_1 \in Z(G)$, a contradiction. If $G \in \mathcal{G}_8$ then $o(z_2) = p^{m_2}$ and $o(z_3) = \infty$. So, $s = q = 0$ and $p|r$ which implies that $d_1 \in \langle t_1 \rangle \times \langle z_2 \rangle \subseteq Z(G)$, a contradiction. If $G \in \mathcal{G}_9$, we have that $o(z_2) = o(z_3) = \infty$, so $r = n = s = q = 0$ which implies that $d_1 \in \langle t_1 \rangle \subseteq Z(G)$, again a contradiction.
\\

Now assume that the second possibility occurs. If $G \in \mathcal{G}_7$ we can suppose that $p\not= 2$ since the case $p=2$ was studied in \cite{casofinito}. Hence,

\[ d_1^p = t_1^{mp}t_2^{r+np}t_3^{s+qp} = \alpha = t_1^at_2^bt_3^c\]
with $o(\alpha) = o(t_1^at_2^bt_3^c) = p^k$. Then $p^{m_1}|(mp-a)$,\; $p^{m_2}|(r+mp-b)$ and \; $p^{m_3}|(s+qp-c)$. Now, $s = t_1^{ajp^{k-1}}t_2^{bjp^{k-1}}t_3^{cjp^{k-1}} = t_1^{p^{m_1-1}}$ which implies that $p^{m_1}|(ajp^{k-1} - p^{m_1-1})$,\; $p^{m_2}|(bjp^{k-1})$ and $p^{m_3}| (cjp^{k-1})$. Since $\langle \alpha \rangle$ is a direct factor of $Z(G)$, it follows that $k = m_1, m_2$ or $m_3$. Note that $p|a$ and $p^{m_1}|(ajp^{k-1} - p^{m_1-1})$ which implies that $k$ cannot be equal to $m_1$.

If $k = m_2$, we have that $p|b$ and since $p^{m_2}|(r+mp-b)$ it follows that $p|r$. So, $r = pr'$ for some positive integer $r'$ and $d_1 = x^ry^st_1^mt_2^nt_3^q = y^st_1^mt_2^{n+r'}t_3^q$. Since $p^{m_3}|(cjp^{k-1})$ we have that either $m_2-1 \geq m_3$ or $p|c$. If $p|c$ then $p|s$ and $d_1 \in Z(G)$, a contradiction. So, $m_2 -1 \geq m_3$. Since $p|b$, we have that $(t_1^at_2^bt_3^c)^{p^{m_2-1}} = t_1^{ap^{m_2-1}}$.

If $m_1 \leq m_2$ then $t_1^{ap^{m_2-1}} = 1$ as $p|a$, a contradiction since $o(t_1^at_2^bt_3^c) = p^{m_2}$. So we can assume that $m_1 > m_2 > m_3$. By Lemma \ref{lema2} there exists $x_1 \in D$ such that $x_1 = xz'$ with $z' = t_1^{A}t_2^{B}t_3^{C} \in Z(G)$. Thus, $x_1^p \in D \cap Z(G) \subseteq Z(D) = \langle \alpha \rangle \times \langle \beta \rangle = \langle t_1^at_2^bt_3^c \rangle \times \langle t_1^{a'}t_2^{b'}t_3^{c'} \rangle$. So,
\[ t_1^{Ap}t_2^{Bp+1}t_3^{Cp} = (t_1^at_2^bt_3^c)^u(t_1^{a'}t_2^{b'}t_3^{'c})^{u'}.\]
Since $p|b$ the equation above implies that $p$ divides neither $b'$ nor $u'$. Since $o(t_1^at_2^bt_3^c) = p^{m_2}$ it follows that $o(t_1^{a'}t_2^{b'}t_3^{c'})= p^{m_1}$ or $p^{m_3}$. First, suppose that $o(\beta) = p^{m_1}$. We have that $p^{m_1} |(au+a'u' - Ap)$ which implies that $p|a'$. Then $(t_1^{a'}t_2^{b'}t_3^{c'})^{p^{m_1-1}} = 1$, a contradiction. If $o(\beta) = p^{m_3}$, then $(t_1^{a'}t_2^{b'}t_3^{c'})^{p^{m_3}} = 1$ which implies that $t_2^{b'p^{m_3}} = 1$, again a contradiction since $o(t_2)= p^{m_2}$, $p$ does not divide $b'$ and $m_3 < m_2$.

If $k = m_3$, then $p|c$ which implies that $p|s$. As for groups in $\mathcal{G}_7$ we have that $y^p = t_3$, writing $s = ps'$ we obtain $d_1 = x^ry^st_1^mt_2^nt_3^q$ = $x^rt_1^mt_2^nt_3^{s'+q}$. Since $p^{m_2}| (bjp^{m_3-1})$ it follows that $p|b$ or $m_3 - 1 \geq m_2$. If $p|b$ we can conclude that $p|r$ which implies that $d_1 \in Z(G)$, a contradiction. So, $m_3 - 1 \geq m_2$. As $p|c$, we have that $(t_1^at_2^bt_3^c)^{p^{m_3-1}} = t_1^{ap^{m_3-1}}$. If $m_3 \geq m_1$ then $t_1^{ap^{m_3-1}} = 1$, again a contradiction, since we are assuming $o(\alpha) = p^{m_3}$. Then $m_2 < m_3 < m_1$. By Lemma \ref{lema2}, there exists $y_1 \in D$ such that $y_1 = yz$ with $z = t_1^{C}t_2^{D}t_3^{E} \in Z(G)$. So, $y_1^p = t_3t_1^{Cp}t_2^{Dp}t_3^{Ep} \in Z(D)$. We can write
\[ t_1^{Cp}t_2^{Dp}t_3^{Ep} = (t_1^at_2^bt_3^c)^v(t_1^{a'}t_2^{b'}t_3^{c'})^{v'}\]

Since $p|c$, it follows that $p$ divides neither $c'$ nor $v'$. Suppose that $o(\beta) = p^{m_1}$. Since $p|a$, we have that $p^{m_1}|(av + a'v' - Cp)$. Hence, $p|a'$ and $(t_1^{a'}t_2^{b'}t_3^{c'})^{p^{m_1-1}} = 1$, which is a contradiction since we are assuming that $o(\beta)= p^{m_1}$. If $o(\beta) = p^{m_2}$, then $(t_1^{a'}t_2^{b'}t_3^{c'})^{p^{m_2}} = 1$ which implies that $t_3^{c'p^{m_2}} = 1$. Since $m_2 < m_3$, it follows that $p| c'$, a contradiction. So, all the groups in $\mathcal{G}_7$ are indecomposable.
\\

If $G \in \mathcal{G}_8$, then $\alpha = t_1^at_2^b$ with $o(t_1^at_2^b)= p^k$. We can write $s = t_1^{ajp^{k-1}}t_2^{bjp^{k-1}}$, where $p$ does not divide $j$. In this case, using equation (\ref{horrible}) we have
\[ t_1^{mp+rs\delta p^{m_1-1}}t_2^{r+np}z_3^{s+qp} = t_1^at_2^b \]
where $\delta = \left\{
                  \begin{array}{ll}
                    0, & \hbox{if $p\not=2$;} \\
                    1, & \hbox{if $p=2$.}
                  \end{array}
                \right.$
\\

Hence, $p^{m_1} |(mp+\delta rsp^{m_1-1} - a)$,\; $p^{m_2}|(r+np-b)$ and $s+qp = 0$. So, $s= q = 0$ and $d_1 = x^rt_1^mt_2^n$. But, $s = t_1^{ajp^{k-1}}t_2^{bjp^{k-1}} = t_1^{p^{m_1-1}}$, which implies $t_2^{bjp^{k-1}} =1$ and then $k -1 \geq m_2$ or $p|b$.
If $p|b$, since $p^{m_2} |(r+np-b)$, it follows that $p|r$ which implies $d_1 \in Z(G)$, a contradiction. So, $k -1 \geq m_2$ which implies $k = m_1$ because $\langle \alpha \rangle $ is a direct factor of order $p^k$ of $Z(G)$, a product of two direct factors of orders $p^{m_i}, i =1,2$ and one factor $\langle z_3 \rangle$ of infinite order. We claim that $a$ is not a multiple of $p$. In fact, if not $(t_1^at_2^b)^{p^{m_1-1}} =1$, which is a contradiction since the order of $\alpha$ is $p^k = p^{m_1}$. Since $p^{m_1}|(mp+\delta rsp^{m_1-1}-a)$ and $m_1 \not= 1$ (as $k = m_1 > m_2$), it follows that $p|a$, a contradiction. Hence, groups in $\mathcal{G}_8$ are indecomposable.
\\

If $G \in \mathcal{G}_9$, it follows that $\alpha = t_1^a$, because
$s \in \langle \alpha \rangle, s^p =1$ and $o(z_2) = o(z_3) = \infty$.
Then $z_2^{r+np} = z_3^{s+qp} = 1$ which implies that $r = n = s = q = 0$
and hence $d_1 = t_1^m \in Z(G)$, a contradiction.
We conclude that groups in $\mathcal{G}_9$ are also indecomposable.



Let $G_i \in \mathcal{G}_i, i = 7,8,9$. Because of cardinality considerations, it readily follows that $G_7 \not\simeq G_8$ and $G_7 \not\simeq G_9$. Since rank of the free part of the center of $G_8$ is 1 and this rank is 2 in $G_9$, we have that $G_8 \not\simeq G_9$. $\hfill{\fimprova}$
\\
\newpage

The results in this section give us the following.
\begin{theorem}
     Let $G$ be a finitely generated indecomposable group such that $G/Z(G) \simeq C_p \times C_p$ where $C_p$ denotes a cyclic group of prime order $p$. Then $G$ is of the form $G = \langle x,y, Z(G) \rangle$ with $x^p, y^p \in Z(G)$ and of one of the nine types of non isomorphic groups listed in table 1. 
\medskip


\begin{center}
\begin{table}[htb]
\begin{center}
\caption{Classification Table}
\begin{tabular}{||c||c||c||c||}
  \hline
  $G$ & $Z(G)$ & $x^p$ & $y^p$ \\
  \hline
  1 & $\langle t_1 \rangle$ & 1 & 1 \\
  2 & $\langle t_1 \rangle$ & $t_1$ & $t_1$ \\
  3 & $\langle t_1 \rangle \times \langle t_2 \rangle$ & 1 & $t_2$ \\
  4 & $\langle t_1 \rangle \times \langle t_2 \rangle$ & $t_1$ & $t_2$ \\
  5 & $\langle t_1 \rangle \times \langle u_1 \rangle$ & 1 & $u_1$ \\
  6 & $\langle t_1 \rangle \times \langle u_1 \rangle$ & $t_1$ & $u_1$ \\
  7 & $\langle t_1 \rangle \times \langle t_2 \rangle \times \langle t_3 \rangle$ & $t_2$& $t_3$ \\
  8 & $\langle t_1 \rangle \times \langle t_2 \rangle \times \langle u_1$  & $t_2$ & $u_1$ \\
  9 & $\langle t_1 \rangle \times \langle u_1 \rangle \times \langle u_2$ & $u_1$ & $u_2$ \\
  \hline
\end{tabular}
\end{center}
\end{table}
\end{center}

In each case of the table 1 we have that $o(t_i) = p^{m_i}, i=1,2,3$ and $o(u_j) = \infty, j = 1,2$.
\end{theorem}



\begin{thebibliography}{}

\bibitem{CG} Chein, O., Goodaire, E.G. (1986) Loops whose loop
rings are alternative, {\slshape Comm. in Algebra} 14:293-310.

\bibitem{good} Goodaire, E.G.(1983). Alternative Loop
Rings, {\slshape Publ. Math. Debrecen 30}, pp. 31-38.


\bibitem{history} Goodaire, E.G. (1999) A Brief History of Loop
Rings, {\slshape Mat. Contemp.} 16: 93-109.

\bibitem{livro} Goodaire,E.G.,Jespers, E., Milies, C.P. (1996)
 Alternative Loop Rings, {\slshape North Holland Mathematics
Studies} 184


\bibitem{casofinito} Jespers, E., Leal, G., Milies, C.P.(1995)
Classifying Indecomposable RA Loops, {\slshape J. Algebra}, 176:569-584.

\bibitem{isomorfismo} Leal, G., Milies, C.P. (1993) Isomorphic Group (and Loop) Algebras, {\slshape J. Algebra}, 155:195-210.



\end{thebibliography}
\end{document}